\documentclass[11pt]{article}
\usepackage{cite}
\usepackage{mathrsfs}
\usepackage{amsfonts}
\usepackage{amsmath}
\usepackage{amsfonts,amssymb}
\usepackage{dsfont}
\usepackage{curves}
\usepackage{mathrsfs}
\usepackage{pifont}
\usepackage{amssymb}
\allowdisplaybreaks

\numberwithin{equation}{section}

\date{}

\begin{document}
\title{Independence number and connectivity for fractional\\ $(a,b,k)$-critical covered graphs
}
\author{\small  Sizhong Zhou$^{1}$\footnote{Corresponding author. E-mail address: zsz\_cumt@163.com (S. Zhou)}, Jiancheng Wu$^{1}$, Hongxia Liu$^{2}$\\
\small $1$. School of  Science, Jiangsu University of Science and Technology\\
\small  Mengxi Road 2, Zhenjiang, Jiangsu 212003, China\\
\small $2$. School of Mathematics and Informational Science, Yantai University\\
\small Yantai, Shandong 264005, P. R. China\\
}

\maketitle
\begin{abstract}
\noindent A graph $G$ is a fractional $(a,b,k)$-critical covered graph if $G-U$ is a fractional $[a,b]$-covered graph for every $U\subseteq V(G)$
with $|U|=k$, which is first defined by Zhou, Xu and Sun (S. Zhou, Y. Xu, Z. Sun, Degree conditions for fractional $(a,b,k)$-critical covered
graphs, Information Processing Letters, DOI: 10.1016/j.ipl.2019.105838). Furthermore, they derived a degree condition for a graph to be a
fractional $(a,b,k)$-critical covered graph. In this paper, we gain an independence number and connectivity condition for a graph
to be a fractional $(a,b,k)$-critical covered graph and verify that $G$ is a fractional $(a,b,k)$-critical covered graph if
$$
\kappa(G)\geq\max\Big\{\frac{2b(a+1)(b+1)+4bk+5}{4b},\frac{(a+1)^{2}\alpha(G)+4bk+5}{4b}\Big\}.
$$
\\
\begin{flushleft}

(2010) Mathematics Subject Classification: 05C70
\end{flushleft}
\end{abstract}

\section{Introduction}

We discuss finite graphs which have neither multiple edges nor loops. Let $G$ be a graph with vertex set $V(G)$ and edge set $E(G)$. Let $X$
and $Y$ be disjoint vertex subsets of $G$. The number of edges of $G$ joining $X$ to $Y$ is denoted by $e_G(X,Y)$. We denote by $G[X]$ and
$G-X$ the subgraph of $G$ induced by $X$ and $V(G)\setminus X$. A vertex subset $X$ of $G$ is called independent if $G[X]$ does not contain
edges. We use $\alpha(G)$ and $\kappa(G)$ to denote the independence number and the connectivity of $G$, respectively. For $x\in V(G)$, we
denote by $d_G(x)$ the degree of $x$ in $G$ and by $N_G(x)$ the set of vertices adjacent to $x$ in $G$. Setting $N_G[x]=N_G(x)\cup\{x\}$.

Let $a\leq b$ be two positive integers. A spanning subgraph $F$ of $G$ is called an $[a,b]$-factor if $a\leq d_F(x)\leq b$ for every $x\in V(G)$.
Let $h:E(G\rightarrow[0,1]$ be a function. If $a\leq\sum\limits_{e\ni x}{h(e)}\leq b$ holds for all $x\in V(G)$, then we call graph $F$ with
vertex set $V(G)$ and edge set $E_h$ a fractional $[a,b]$-factor of $G$ with indicator function $h$, where $E_h=\{e:e\in E(G),h(e)>0\}$. A
fractional $[a,b]$-factor of $G$ is called a fractional $r$-factor of $G$ if $a=b=r$.  A graph $G$ is called a fractional $[a,b]$-covered
graph if $G$ admits a fractional $[a,b]$-factor containing any given edge $e$. Bian and Zhou \cite{BZ} posed an independence number
and connectivity condition for the existence of fractional $r$-factors in graphs, which is a special case of Bian and Zhou's result.

\medskip

\noindent{\textbf{Theorem 1}} (\cite{BZ}). Let $G$ be a graph, and $r\geq1$ be an integer. Then $G$ possesses a fractional $r$-factor if
$$
\kappa(G)\geq\max\Big\{\frac{(r+1)^{2}}{2},\frac{(r+1)^{2}\alpha(G)}{4r}\Big\}.
$$

\medskip

Some other results related to factors \cite{BY,KO,KL,Zi,ZSX,Zs,Zr,ZYX} and fractional factors \cite{GGC,LZc,K,ZS,ZSo,ZXX,Za,ZZ,YHa} of
graphs were obtained by many authors. A graph $G$ is called a fractional $(a,b,k)$-critical covered graph if after removing any $k$
vertices of $G$, the resulting graph of $G$ is a fractional $[a,b]$-covered graph, which is first defined by Zhou, Xu and Sun \cite{ZXS}.
Furthermore, Zhou, Xu and Sun \cite{ZXS} obtained a degree condition for graphs being fractional $(a,b,k)$-critical covered.

\medskip

\noindent{\textbf{Theorem 2}} (\cite{ZXS}). Let $a,b$ and $k$ be three integers with $a\geq1$, $b\geq\max\{2,a\}$ and $k\geq0$, and let $G$
be a graph of order $n$ with $n\geq\frac{(a+b)(a+b-1)+bk+3}{b}$ and $\delta(G)\geq a+k+1$. Then $G$ is a fractional $(a,b,k)$-critical
covered graph if
$$
\max\{d_G(x),d_G(y)\}\geq\frac{an+bk+2}{a+b}
$$
for every pair of nonadjacent vertices $x$ and $y$ of $G$.

\medskip

In this paper, we study the relationship between independence number, connectivity and fractional $(a,b,k)$-critical covered graphs, and
gain a new result on the existence of fractional $(a,b,k)$-critical covered graphs.

\medskip

\noindent{\textbf{Theorem 3.}} Let $a,b$ and $k$ be three integers with $b\geq a\geq1$ and $k\geq0$, and let $G$ be a graph. Then $G$ is a
fractional $(a,b,k)$-critical covered graph if
$$
\kappa(G)\geq\max\Big\{\frac{2b(a+1)(b+1)+4bk+5}{4b},\frac{(a+1)^{2}\alpha(G)+4bk+5}{4b}\Big\}.
$$

\medskip

We immediately derive the following result from Theorem 3.

\medskip

\noindent{\textbf{Corollary 1.}} Let $a$ and $b$ be integers with $b\geq a\geq1$, and let $G$ be a graph. Then $G$ is a fractional $[a,b]$-covered
graph if
$$
\kappa(G)\geq\max\Big\{\frac{2b(a+1)(b+1)+5}{4b},\frac{(a+1)^{2}\alpha(G)+5}{4b}\Big\}.
$$

\section{Proof of Theorem 3}

We use the following lemma to verify Theorem 3.

\medskip

\noindent{\textbf{Lemma 1}} (\cite{LYZ}). Let $G$ be a graph, and let $a,b$ be two integers with $b\geq a\geq0$. Then $G$ is a fractional
$[a,b]$-covered graph if and only if for any vertex subset $X$ of $G$,
$$
\theta_G(X,Y)=b|X|+d_{G-X}(Y)-a|Y|\geq\varepsilon(X,Y)
$$
where $Y=\{y:y\in V(G)\setminus X, d_{G-X}(y)\leq a\}$ and $\varepsilon(X,Y)$ is defined by
\[
 \varepsilon(X,Y)=\left\{
\begin{array}{ll}
2,&if \ X \ is \ not \ independent,\\
1,&if \ X \ is \ independent \ and \ there \ is \ an \ edge \ joining \ X \ and \
V(G)\setminus(X\cup Y), \ or\\
&there \ is \ an \ edge \ e=xy \ joining \ X \ and \ Y \ with \ d_{G-X}(y)=a \ for \ y\in Y,\\
0,&otherwise.\\
\end{array}
\right.
\]

\medskip

\noindent{\it Proof of Theorem 3.} \ Let $U\subseteq V(G)$ with $|U|=k$, and let $H=G-U$. It suffices to verify that $H$ is a fractional
$[a,b]$-covered graph. We shall prove this by contradiction. Suppose that $H$ is not a fractional $[a,b]$-covered graph. Then using Lemma
1, we have
$$
\theta_H(X,Y)=b|X|+d_{H-X}(Y)-a|Y|\leq\varepsilon(X,Y)-1\eqno(1)
$$
for some vertex subset $X$ of $H$ and $Y=\{y:y\in V(H)\setminus X:d_{H-X}(y)\leq a\}$.

\noindent{\bf Claim 1.} \ $Y\neq\emptyset$.

\noindent{\it Proof.} \ Let $Y=\emptyset$. Then by (1) and $\varepsilon(X,\emptyset)\leq|X|$, we get
$$
\varepsilon(X,\emptyset)-1\geq\theta_H(X,\emptyset)=b|X|\geq|X|\geq\varepsilon(X,\emptyset),
$$
this is a contradiction. Claim 1 is proved. \hfill $\Box$

\noindent{\bf Claim 2.} \ $X\neq\emptyset$.

\noindent{\it Proof.} \ Let $X=\emptyset$. Then $\varepsilon(\emptyset,Y)=0$. Note that $\delta(H)=\delta(G-U)\geq\delta(G)-k\geq\kappa(G)-k
\geq\frac{2b(a+1)(b+1)+4bk+5}{4b}-k=\frac{(a+1)(b+1)}{2}+\frac{5}{4b}>a+1$. Then using (1) and Claim 1 we obtain
\begin{eqnarray*}
-1&\geq&\theta_H(\emptyset,Y)=d_H(Y)-a|Y|\geq(\delta(H)-a)|Y|\\
&>&((a+1)-a)|Y|=|Y|\geq1,
\end{eqnarray*}
which is a contradiction. We prove Claim 2. \hfill $\Box$

Since $Y\neq\emptyset$ (by Claim 1), we take $y_1\in Y$ such that $y_1$ is the vertex with the minimum degree in $H[Y]$. Setting $N_1=N_H[y_1]\cap Y$
and $Y_1=Y$. If $Y-\bigcup\limits_{1\leq j<i}N_j\neq\emptyset$ for $i\geq2$, write $Y_i=Y-\bigcup\limits_{1\leq j<i}N_j$. Then take $y_i\in Y_i$ such
that $y_i$ is the vertex with the minimum degree in $H[Y_i]$, and $N_i=N_H[y_i]\cap Y_i$. We go on these procedures until we arrive at the situation
in which $Y_i=\emptyset$ for some $i$, say for $i=m+1$. Then from the definition above we see that $\{y_1,y_2,\cdots,y_m\}$ is an independent set of
$H$, and is also an independent set of $G$. Obviously, $m\geq1$ by Claim 1.

Write $|N_i|=n_i$, then $|Y|=\sum\limits_{1\leq i\leq m}n_i$. Setting $W=V(H)\setminus(X\cup Y)$ and $\kappa(H-X)=t$.

\noindent{\bf Claim 3.} \ $m\neq1$ or $W\neq\emptyset$.

\noindent{\it Proof.} \ Let $m=1$ and $W=\emptyset$. Then we easily see that
$$
|V(G)|=|U|+|X|+n_1=k+|X|+n_1>\kappa(G)\geq\frac{2b(a+1)(b+1)+4bk+5}{4b},
$$
namely,
$$
|X|>\frac{2b(a+1)(b+1)+5}{4b}-n_1.\eqno(2)
$$

According to (1), (2), $\varepsilon(X,Y)\leq2$ and the choice of $y_1$, we get
\begin{eqnarray*}
1&\geq&\varepsilon(X,Y)-1\geq\theta_H(X,Y)=b|X|+d_{H-X}(Y)-a|Y|\\
&=&b|X|+n_1(n_1-1)-an_1\\
&>&b\Big(\frac{2b(a+1)(b+1)+5}{4b}-n_1\Big)+n_1(n_1-1)-an_1\\
&=&\frac{2b(a+1)(b+1)+5}{4}+\Big(n_1-\frac{a+b+1}{2}\Big)^{2}-\frac{(a+b+1)^{2}}{4}\\
&\geq&\frac{2b(a+1)(b+1)+5}{4}-\frac{(a+b+1)^{2}}{4}\\
&=&\frac{b^{2}(2a+1)-(a+1)^{2}+5}{4}\\
&\geq&\frac{a^{2}(2a+1)-(a+1)^{2}+5}{4}\\
&=&\frac{a(a^{2}-1)+2}{2}\geq1,
\end{eqnarray*}
which is a contradiction. We verify Claim 3. \hfill $\Box$

\noindent{\bf Claim 4.} \ $d_{H-X}(Y)\geq\sum\limits_{1\leq i\leq m}n_i(n_i-1)+\frac{mt}{2}$.

\noindent{\it Proof.} \ In terms of the choice of $y_i$, we derive
$$
\sum_{1\leq i\leq m}(\sum_{y\in N_i}d_{Y_i}(y))\geq\sum_{1\leq i\leq m}n_i(n_i-1).\eqno(3)
$$
For the left-hand side of (3), an edge joining a vertex $x$ in $N_i$ and a vertex $y$ in $N_j$ ($i<j$) is counted only once, namely, it is
counted in $d_{Y_i}(x)$ but not in $d_{Y_j}(y)$. Hence, we admit
$$
d_{H-X}(Y)\geq\sum_{1\leq i\leq m}n_i(n_i-1)+\sum_{1\leq i<j\leq m}e_H(N_i,N_j)+e_H(Y,W).\eqno(4)
$$

In light of $\kappa(H-X)=t$ and Claim 3, we get
$$
e_H(N_i,\bigcup_{j\neq i}N_j)+e_H(N_i,W)\geq t\eqno(5)
$$
for every $N_i$ ($1\leq i\leq m$). Using (5), we obtain
$$
\sum_{1\leq i\leq m}(e_H(N_i,\bigcup_{j\neq i}N_j)+e_H(N_i,W))=2\sum_{1\leq i<j\leq m}e_H(N_i,N_j)+e_H(Y,W)\geq mt.
$$
We easily see that
$$
\sum_{1\leq i<j\leq m}e_H(N_i,N_j)+e_H(Y,W)\geq\frac{mt}{2}.\eqno(6)
$$

It follows from (4) and (6) that
$$
d_{H-X}(Y)\geq\sum_{1\leq i\leq m}n_i(n_i-1)+\frac{mt}{2}.
$$
We finish the proof of Claim 4. \hfill $\Box$\

Note that $|Y|=\sum\limits_{1\leq i\leq m}n_i$. It follows from (1), $\varepsilon(X,Y)\leq2$ and Claim 4 that
\begin{eqnarray*}
1&\geq&\varepsilon(X,Y)-1\geq\theta_H(X,Y)=b|X|+d_{H-X}(Y)-a|Y|\\
&\geq&b|X|+\sum\limits_{1\leq i\leq m}n_i(n_i-1)+\frac{mt}{2}-a\sum\limits_{1\leq i\leq m}n_i\\
&=&b|X|+\sum\limits_{1\leq i\leq m}\Big(\Big(n_i-\frac{a+1}{2}\Big)^{2}-\frac{(a+1)^{2}}{4}\Big)+\frac{mt}{2}\\
&\geq&b|X|-\frac{m(a+1)^{2}}{4}+\frac{mt}{2}\\
&=&b|X|+\Big(\frac{t}{2}-\frac{(a+1)^{2}}{4}\Big)m,
\end{eqnarray*}
namely,
$$
1\geq b|X|+\Big(\frac{t}{2}-\frac{(a+1)^{2}}{4}\Big)m.\eqno(7)
$$

Combining (7), Claim 2 and $m\geq1$, we admit
\begin{eqnarray*}
1&\geq&b|X|+\Big(\frac{t}{2}-\frac{(a+1)^{2}}{4}\Big)m\\
&\geq&b+\Big(\frac{t}{2}-\frac{(a+1)^{2}}{4}\Big)m\\
&\geq&1+\Big(\frac{t}{2}-\frac{(a+1)^{2}}{4}\Big)m,
\end{eqnarray*}
which implies
$$
\frac{t}{2}-\frac{(a+1)^{2}}{4}\leq0.\eqno(8)
$$

We easily see that $\alpha(G)\geq\alpha(G[Y])=\alpha(H[Y])\geq m$ and $\kappa(G)\leq\kappa(G-U)+k=\kappa(H)+k\leq\kappa(H-X)+|X|+k=t+|X|+k$.
In light of $\kappa(G)\geq\max\Big\{\frac{2b(a+1)(b+1)+4bk+5}{4b},\frac{(a+1)^{2}\alpha(G)+4bk+5}{4b}\Big\}$, (7) and (8), we derive
\begin{eqnarray*}
1&\geq&b|X|+\Big(\frac{t}{2}-\frac{(a+1)^{2}}{4}\Big)m\\
&\geq&b(\kappa(G)-k-t)+\Big(\frac{t}{2}-\frac{(a+1)^{2}}{4}\Big)\alpha(G)\\
&\geq&b(\kappa(G)-k-t)+\Big(\frac{t}{2}-\frac{(a+1)^{2}}{4}\Big)\cdot\frac{4b\kappa(G)-4bk-5}{(a+1)^{2}}\\
&=&\frac{5}{4}+\Big(\frac{4b\kappa(G)-4bk-5}{2(a+1)^{2}}-b\Big)t\\
&\geq&\frac{5}{4}+\Big(\frac{b(a+1)(b+1)}{(a+1)^{2}}-b\Big)t\\
&\geq&\frac{5}{4},
\end{eqnarray*}
it is a contradiction. Theorem 3 is verified. \hfill $\Box$\\\

\medskip

\section{Remarks}

\noindent{\bf Remark 1.} \ Now we discuss a sharpness of the connectivity condition in Theorem 3. This condition is best possible in the sense
which we cannot replace by $\kappa(G)\geq\frac{2b(a+1)(b+1)+4bk+4}{4b}$.

Let $a=b=1$ and $k\geq0$ be an integer. We construct a graph $G=K_{3+k}\vee(2K_1)$. Then it is obvious that $\kappa(G)=3+k=\frac{2b(a+1)(b+1)+4bk+4}{4b}$
by the definition of $\kappa(G)$. We select $U\subseteq V(K_{3+k})$ with $|U|=k$, $X=V(K_{3+k})\setminus U$, $Y=V(2K_1)$ and $H=G-U$. Clearly,
$\varepsilon(X,Y)=2$. Thus, we derive
$$
\theta_H(X,Y)=b|X|+d_{H-X}(Y)-a|Y|=3-2=1<2=\varepsilon(X,Y).
$$
By Lemma 1, $H$ is not a fractional $[a,b]$-covered graph, namely, $G$ is not a fractional $(a,b,k)$-critical covered graph.

\medskip

\noindent{\bf Remark 2.} \ The condition $\kappa(G)\geq\frac{(a+1)^{2}\alpha(G)+4bk+5}{4b}$ in Theorem 3 is best possible.

Let $b\geq a\geq1$, $m\geq1$ and $k\geq0$ be integers such that $a$ is odd and $\frac{(a+1)^{2}m+4}{4b}$ is an integer. We construct a graph
$G=K_{p+k}\vee(mK_\frac{a+1}{2})$, where $p=\frac{(a+1)^{2}m+4}{4b}$. From the definitions of $\alpha(G)$ and $\kappa(G)$, we easily see that
$\alpha(G)=m$ and $\kappa(G)=p+k=\frac{(a+1)^{2}m+4bk+4}{4b}=\frac{(a+1)^{2}\alpha(G)+4bk+4}{4b}$. We select $U\subseteq V(K_{p+k})$ with
$|U|=k$, $X=V(K_{p+k})\setminus U$, $Y=V(mK_\frac{a+1}{2})$ and $H=G-U$. Obviously, $\varepsilon(X,Y)=2$. Thus, we gain
\begin{eqnarray*}
\theta_H(X,Y)&=&b|X|+d_{H-X}(Y)-a|Y|\\
&=&b\cdot\frac{(a+1)^{2}m+4}{4b}+m\cdot\frac{a+1}{2}\cdot\Big(\frac{a+1}{2}-1\Big)-am\cdot\frac{a+1}{2}\\
&=&1<2=\varepsilon(X,Y).
\end{eqnarray*}
Using Lemma 1, $H$ is not a fractional $[a,b]$-covered graph, and so $G$ is not a fractional $(a,b,k)$-critical covered graph.


\end{document}